\documentclass{amsart}

\usepackage{amsmath,amsfonts,amssymb,amsthm,amscd,latexsym,euscript,bbm}

   \newtheorem{lemma}{Lemma}[section]  
   \newtheorem{kor}[lemma]{Corollary}        
   \newtheorem{satz}[lemma]{Theorem}         

\theoremstyle{definition}
\newtheorem{Def}[lemma]{Definition}                                          

    \newtheorem{Bem}[lemma]{Remark}        

\usepackage[all]{xy}

\newenvironment{beweis}{\noindent\textbf{Proof:}}{\par \hfill $\Box$ \hspace{1cm} \par}

\numberwithin{equation}{section}        
\newcommand{\Ref}[1]{(\ref{#1})}        

\newcounter{zahl}%
\newenvironment{punkt}{\begin{list}{{\rm{(\roman{zahl})}}}%
    {\usecounter{zahl}%
     \setlength{\leftmargin}{0pt} \setlength{\itemindent}{4pt} \setlength{\topsep}{2pt} \setlength{\parsep}{2pt} }}%
    {\end{list}}%

\newcommand{\co}{\colon\thinspace}
\newcommand{\mc}[1]{\ensuremath{\mathcal{#1}}}
\newcommand{\aus}{\raisebox{1pt}{\ensuremath{\,{\scriptstyle\in}\,}}}%

\newcommand{\toh}[1]{\ensuremath{\stackrel{#1}{\rightarrow}}}%
\newcommand{\colim}{\operatorname*{colim}}
\newcommand{\cosk}{\operatorname{cosk}}
\newcommand{\frei}{\,\_\!\_\,}
\newcommand{\id}{\operatorname{id}}
\newcommand{\inj}{\ensuremath{\hookrightarrow}}

\newcommand{\hrl}{\ensuremath{\raisebox{1pt}{(}}}
\newcommand{\hrr}{\ensuremath{\raisebox{1pt}{)}}}
\newcommand{\ho}[1]{\ensuremath{H\!o({#1})}}

\newcommand{\bth}{\raisebox{1pt}{\ensuremath{\,\scriptstyle \geq\,}}}%
\newcommand{\sth}{\raisebox{1pt}{\ensuremath{\,\scriptstyle \leq\,}}}%
\newcommand{\uber}{\raisebox{1pt}{\ensuremath{\,{\scriptstyle > }\,}}}%
\newcommand{\unter}{\raisebox{1pt}{\ensuremath{\,{\scriptstyle < }\,}}}%
\newcommand{\ent}{\raisebox{1pt}{\ensuremath{\,{\scriptstyle \subseteq }\,}}}%
\newcommand{\pullback}{\displaystyle \cdot\!\!\lrcorner}%
\newcommand{\ol}[1]{\overline{#1}}
\newcommand{\wt}[1]{\widetilde{#1}}
\newcommand{\Hom}[3]{\ensuremath{{\rm Hom}_{#1}\hrl#2,#3\hrr}}%
\newcommand{\map}{\ensuremath{{\rm map}}}
\renewcommand{\hom}{\ensuremath{{\rm hom}}}

\newcommand{\diag}[2]{ \begin{align} \begin{split} \xymatrix{#1} \end{split} \label{#2} \end{align}}%

\newcommand{\diagr}[1]{ \begin{equation*} \xymatrix{#1} \end{equation*}}%

\newcommand{\pre}{\ensuremath{s{\rm Pre}_{\mc{C}}}}%
\newcommand{\shv}{\ensuremath{s{\rm Shv}_{\mc{C}}}}%
\newcommand{\pregd}{\ensuremath{s{\rm Pre}_{\mc{C}}({\rm Gd})}}%
\newcommand{\shvgd}{\ensuremath{s{\rm Shv}_{\mc{C}}({\rm Gd})}}%
\begin{document}

\title  [On the homotopy theory of $n$-types]
        {On the homotopy theory of $n$-types}
\author{Georg Biedermann}

\address{Department of Mathematics, Middlesex College, The University of Western Ontario, London, Ontario N6A 5B7, Canada}

\email{gbiederm@uwo.ca}



\date{\today}
\dedicatory{}
\commby{}

\begin{abstract}
An $n$-truncated model structure on simplicial (pre-)sheaves is described having as weak equivalences maps that induce isomorphisms on certain homotopy sheaves only up to degree $n$. Starting from one of Jardine's intermediate model structures we construct such an $n$-type model structure via Bousfield-Friedlander localization and exhibit useful generating sets of trivial cofibrations. Injectively fibrant objects in these categories are called $n$-hyper\-stacks. The whole setup can consequently be viewed as a description of the homotopy theory of higher hyperstacks. More importantly, we construct analogous $n$-truncations on simplicial groupoids and prove a Quillen equivalence between these settings. We achieve a classification of $n$-types of simplicial presheaves in terms of $(n-1)$-types of presheaves of simplicial groupoids. 
Our classification holds for general $n$. Therefore this can also be viewed as the homotopy theory of (pre-)sheaves of (weak) higher groupoids. 
\end{abstract}

\maketitle
\tableofcontents

\section{Introduction}

A homotopy $n$-type is the homotopy type of a space $X$ whose homotopy groups $\pi_sX$ of degree $s\uber n$ vanish. ``Space'' in our context means simplicial set or, more generally, (pre-)sheaf of simplicial sets on a small Grothendieck site \mc{C}. A good way to describe homotopy theory is via Quillen model structures, see \cite{DSpa:model}. The right notion of equivalence of this structure is a Quillen equivalence.
The goal of this article is to give sense to the following diagram of Quillen equivalences:
\diag{ \pre^n \ar@<2pt>[r]^-{G}\ar@<2pt>[d]^{L^2} & \pregd^{n-1} \ar@<2pt>[l]^-{\ol{W}}\ar@<2pt>[d]^{L^2} \\
       \shv^n \ar@<2pt>[r]^-{G}\ar@<2pt>[u]^{i} & \shvgd^{n-1} \ar@<2pt>[l]^-{\ol{W}}\ar@<2pt>[u]^{i}}{main}
In this diagram \pre\, denotes the category of simplicial presheaves (aka. presheaves of simplicial sets) on the site \mc{C}, \shv\, denotes the category of simplicial sheaves, \pregd\, the category of presheaves of groupoids enriched in simplicial sets and \shvgd\, is the category of sheaves of such groupoids. 
The functor $L^2$ is the sheafification functor, $i$ is the forgetful functor, $G$ is the loop groupoid functor and $\ol{W}$ is the classifying space or universal cocycle functor. These last two functors were introduced in \cite{DK:simp-groupoids}, see also \cite[V.7.]{GoJar:simp}.
The $n$ in this diagram refers to a truncated model structure, where the weak equivalences only take homotopy groups up to degree $n$ into account. 

For $\mc{C}=\ast$ and $n=\infty$ this reduces to the well known Quillen equivalence between simplicial sets and simplicial groupoids proved in \cite{DK:simp-groupoids}. But even for $\mc{C}=\ast$ and $n\unter\infty$ on one side or $n=\infty$ and \mc{C} arbitrary on the other side diagram \Ref{main} seems to be new.

The homotopy theory of simplicial (pre-)sheaves is well prepared ground, see remark \ref{all serve well} for some references. There are several model structures including the injective or global structure and the projective or local structure, all Quillen equivalent to each other. We will see, that we can start from any of these so-called intermediate model structures and obtain a diagram \Ref{main}, all of them Quillen equivalent to each other. So there is plenty of freedom.

The globally or injectively fibrant objects in the $n$-truncated model structure on \pre\, or \shv\, are actually sectionwise $n$-types. 
Hence, for finite $n$ the theory amounts to a classification of (pre-)sheaves of $n$-types in terms of (pre-)sheaves of truncated groupoids enriched in simplicial sets. 
These injectively fibrant models of $n$-types of simplicial (pre-)sheaves have received some attention in recent publications as $n$-(hyper-)stacks, see \cite{Lurie:topoi} and \cite{TV:hag1}. 
So the above diagram gives different ways of describing the homotopy theory of higher hyperstacks, of which the right hand side seems to be completely new. Also on the left hand side this technique offers new insights, e.g. we identify sets of generating trivial cofibrations relative to the underlying model structures on \pre\, or \shv.

I begin in section \ref{section:BF} with a very brief overview of a localization process from \cite{BF:gamma} and improved on in \cite{Bou:telescopic}, which I will use later on to obtain the truncations.

Then in section \ref{section:presheaves}, I quickly recast the necessary homotopy theory on \pre\, and \shv.
There are several choices for such underlying model structures exhibited in \cite{jar:intermediate} called intermediate model structures, because they sit in between the well known projective and injective model structures. Each of these structures serves equally well as a starting point for our theory. 
The corresponding theories on the groupoid side are developed in section \ref{section:groupoids}. It is an important and quite intricate point to not neglect set-theoretic difficulties. 
We will not be concerned with it though, since Jardine in \cite{jar:simplicial-presheaves} has dealt with this problem by considering the concept of small sites and I simply adopt his point of view. See also \cite{jar:lec-simp-pre}.

In section \ref{section:Postnikov-Lokalisierung} the Bousfield-Friedlander localization from section \ref{section:BF} is employed to obtain a model structure where the weak equivalences are given by maps that induce isomorphisms on all homotopy sheaves up to degree $n$ and where the fibrant objects are exactly the simplicial presheaves whose homotopy groups above degree $n$ vanish. Fibrant approximation is given by taking $n$-th Postnikov stages. 
The existence of such a model structure was certainly folklore for a long time and it has been described in \cite{Hir:loc} or in \cite{Lurie:topoi} and \cite{TV:hag1} as the left Bousfield localization along the map $\partial\Delta^{n+2}\to\Delta^{n+2}$. In this article the structure is obtained by Bousfield-Friedlander localization, \cite{BF:gamma} and \cite{Bou:telescopic}, along the $n$-th Postnikov section functor $P_n$. What we gain -- apart from a much simpler construction -- is an explicit description of the class of fibrations. From that we obtain a set of generating trivial cofibrations which is very small. 

One also obtains an analogous truncation on the side of groupoids. Then it is easy to establish the diagram of Quillen equivalences on the first page. All of this is explained in sections \ref{section:groupoids} and \ref{section:truncated groupoids}. 
I argue that the truncated model structures on groupoids enriched in simplicial sets are a good way of avoiding horrendous technical difficulties with higher groupoids. In fact, from this point of view the theory described here is the homotopy theory of (pre-)sheaves of higher groupoids.
\\

I would like to thank Rick Jardine for providing so much input via lectures, lecture notes, preprints and conversations.
Being in London, Ontario, I simply could not help starting to think about simplicial presheaves.

\section{Bousfield-Friedlander localization}
\label{section:BF}

Bousfield-Friedlander localization is a setup to localize a right proper model category along a coaugmented functor satisfying some reasonable axioms. It was used in \cite{BF:gamma} to introduce a model structure on spectra, which has the stable homotopy category as its homotopy category. The axioms there seemed complicated and somehow tailored for the previous use, but the situation was greatly improved in \cite{Bou:telescopic} which makes this way of localizing model categories extremely efficient. It has the advantage of being simple to prove and the resulting model structure remains under good control: simplicial enrichments and properness properties are inherited and fibrations can be characterized quite explicitly.

Let \mc{M} be a right proper model category. Let $Q\co\mc{M}\to\mc{M}$ be a functor with coaugmentation $\eta\co\id\to Q$. Before we give the axioms, that $Q$ and $\eta$ have to satisfy, we describe what is supposed to become the localized model structure.

\begin{Def}
A map $f\co X\to Y$ in \mc{M} is called 
\begin{punkt}
    \item
a {\bf\boldmath $Q$-equivalence} if $Qf\co QX\to QY$ is a weak equivalence.
    \item
a {\bf\boldmath$Q$-fibration} if $f$ has the right lifting property with respect to cofibrations that are also $Q$-equivalences.
\end{punkt}
The new cofibrations are simply the old ones. We call these classes of maps the {\bf\boldmath$Q$-structure} on \mc{M} and denote the structure by $\mc{M}^Q$.
\end{Def}

\begin{Def} \label{BF-axioms}
We will refer several times to the following axioms: 
\begin{punkt}
    \item[{\bf(A.4)}]
$Q$ preserves weak equivalences.
    \item[{\bf(A.5)}]
The maps $\eta(QX)$ and $Q(\eta x)\co QQ(X)\co QX\to QQX$ are weak equivalences.
    \item[{\bf(A.6)}]
Consider the following diagram:
   $$\xy
   \xymatrix"*"@=13pt{ QA\ar[rr] \ar'[d] [dd] & & QX \ar[dd] \\
                                  & &                            \\
              QB \ar'[r]^-<<{\simeq} [rr] & & QY  }    
   \POS(-10,-8)
   \xymatrix@=17pt{ A \ar[dd]_{g} \ar["*"] \ar[rr]^-<<<<{\alpha} & & X \ar@{->>}[dd]^-<<<{f} \ar["*"]_-\simeq\\
                                  & &                            \\
                    B \ar[rr]_{\beta} \ar["*"] & &  Y  \ar["*"]_-\simeq\ar@{}[uull]|->>{\pullback}  }
     \endxy
   $$
Here $f$ is a fibration between fibrant objects, such that $\eta(X)\co X\to QX$ and $\eta(Y)\co Y\to QY$ are weak equivalences. The front square is a pullback and $QB\to QY$ is a weak equivalence. Then $QA\to QX$ is a weak equivalence.
\end{punkt}
\end{Def}

The last axiom {\bf (A.6)} simply asserts that the localized model structure $\mc{M}^Q$ will be right proper. 
The fact that we can assume $\eta(X)$ and $\eta(Y)$ to be equivalences is one of the improvements in \cite[9.4.]{Bou:telescopic}.
We only have to test right properness on local objects and an object $X$ is local if and only if it is fibrant and $\eta(X)$ is a weak equivalence. Another far more important improvement was the removal of a dual axiom of {\bf (A.6)}. The following theorem is taken from \cite[9.3.]{Bou:telescopic}.

\begin{satz} \label{Bousfield-Friedlander}
For a right proper model category \mc{M} with a coaugmented functor $Q$ that satisfies {\bf (A.4)}, {\bf (A.5)} and {\bf (A.6)} the $Q$-structure is a right proper model structure on \mc{M}. If the original model structure on \mc{M} is left proper or simplicial, so is the $Q$-structure. If \mc{M} has functorial factorization, so has $\mc{M}^Q$. A map $p\co X\to Y$ is a $Q$-fibration if and only if it is a fibration such that the following square
\diagr{ X\ar[r]^-{\eta(X)}\ar[d]_f & QX \ar[d]^{Qf} \\
        Y \ar[r]_-{\eta(Y)} & QY }
is a homotopy pullback square.
\end{satz}

\section{Simplicial presheaves and sheaves}
\label{section:presheaves}

In this section we are going to review very briefly the theory of simplicial presheaves and sheaves as outlined in \cite{jar:lec-simp-pre}.

Let \mc{C} be a small Grothendieck site. A presheaf with values in a category \mc{D} is a contravariant functor from \mc{C} to \mc{D}. For the category of such presheaves we write ${\rm Pre}_{\mc{C}}(\mc{D})$. 
Let \mc{S} denote the category of simplicial sets.
A simplicial presheaf $X$ is a functor $\mc{C}^{\rm op}\to\mc{S}$ and we denote the respective category by $s{\rm Pre}_{\mc{C}}:={\rm Pre}_{\mc{C}}(\mc{S})$. 
We denote the category of simplicial sheaves with values in \mc{D} by $s{\rm Shv}_{\mc{C}}(\mc{D})$ and let $s{\rm Shv}_{\mc{C}}=s{\rm Shv}_{\mc{C}}(\mc{S})$.
These categories support various model structures. But first of all we observe that we have simplicial enrichments.

\begin{Bem} \label{simplicial structure}
The categories \pre\, and \shv\, are simplicially enriched, tensored and cotensored with the following definitions: For $X$ and $Y$ in \pre\, or \shv, $K$ in \mc{S} and $U\aus\mc{C}$ we set:
\begin{align*}
    (X\otimes K)(U) & := X(U)\times K \\
    \map(X,Y)_n     & := \Hom{\pre}{X\otimes\Delta^n}{Y} \\
    \hom(K,X)(U)    & := \map_{\mc{S}}(K,X(U))
\end{align*}
In the sheaf case the tensor is the associated sheaf to the presheaf described above.
All the model structure we will consider on \pre\, or \shv\, are simplicial model structures when equipped with the above functors.
\end{Bem}

Let us first concentrate on the presheaves.
To describe the weak equivalences we have to introduce the right analogue of homotopy groups. We define homotopy groups of simplicial sets by
       $$ \pi_s(X,x) := \pi_s({\rm Ex}^{\infty}X,x) $$
where on the right side we take the simplicial homotopy group of the fibrant replacement ${\rm Ex}^{\infty}X$ of $X$. 
Obviously, the ${\rm Ex}^{\infty}$-functor can be promoted to an objectwise fibrant replacement functor on the category \pre.

\begin{Def}
To a simplicial presheaf $X\co\mc{C}^{\rm op}\to\mc{S}$ we can associate two presheaves of sets, the presheaf of vertices $X_0$ and the presheaf of components $\pi_0X$. 
Further we have for every $U\aus\mc{C}$, $s\bth 1$ and $x\aus X_0(U)$ a group $\pi_s(X(U),x)$, which we can assemble into a group object of presheaves over $X_0$ in the following way: Let
    $$ \pi_sX(U) := \coprod_{x\in X_0(U)}\pi_s({\rm Ex}^{\infty}X(U),x),$$
where the projection to $X_0(U)$ takes the summand $\pi_s(X(U),x)$ to $x$. For $s\bth 1$ these gadgets are group objects fibred over $X_0$, which are abelian for $s\bth 2$. We denote the associated group object of sheaves fibred over $X_0$ by $\wt{\pi}_sX$ and call it the {\bf\boldmath$s$-th homotopy sheaf} of $X$, where we often leave the map to $X_0$ to be understood.
\end{Def}

\begin{Def} \label{local weak equivalences}
A map $f\co X\to Y$ is called a {\bf local weak equivalences} if it induces isomorphisms
    $$ \wt{\pi}_0X \to \wt{\pi}_0Y $$
and pullback diagrams of sheaves
\diagr{ \wt{\pi}_sX \ar[r]\ar[d] & \wt{\pi}_sY \ar[d] \\ \wt{X}_0 \ar[r] & \wt{Y}_0}
for all $s\bth 1$. 
This last condition can also be rephrased by saying that for $s\bth 1$ the canonical map
    $$ \pi_sX = \coprod_{x\in X_0}\pi_s(X,x)\to\coprod_{x\in X_0}\pi_s(Y,fx)=:\pi Y_f $$
induces an isomorphism of associated sheaves fibred over $X_0$. This corresponds for the case $\mc{C}=*$ to the usual requirement, that $X\to Y$ induces an isomorphism
    $$ \pi_s(X,x)\to\pi_s(Y,fx)$$
for all $s\bth 0$ and all basepoints $x\aus X$. The vanishing of a homotopy group $\pi_s(X,x)$ for all basepoints $x$ in the classical setting $\mc{C}=*$ corresponds to the fact that the canonical map
    $$ \wt{\pi}_s(X)\to\wt{X}_0 $$
is an isomorphism of sheaves.
\end{Def}

All the model structures we will consider on \pre\, will have local weak equivalences as their equivalences. In particular, they will all have equivalent homotopy categories associated to them. The difference lies in the choice of fibrations and cofibrations. 

\begin{Def} \label{Def. of cofibrations}
\begin{punkt}
   \item
A map $X\to Y$ is called a {\bf projective} or {\bf local} if it induces a fibration $X(U)\to Y(U)$ of simplicial sets in each section $U\aus\mc{C}$.
   \item
A map $X\to Y$ is called a {\bf projective} or {\bf local cofibration} if it has the left lifting property with respect to all objectwise fibrations that are also local weak equivalences.
   \item
A map $X\to Y$ is called an {\bf injective} or {\bf global} if it induces a cofibration $X(U)\to Y(U)$ of simplicial sets in each section $U\aus\mc{C}$.
   \item
A map $X\to Y$ is called an {\bf injective} or {\bf global fibration} if it has the right lifting property with respect to all injective cofibrations that are also local weak equivalences.
\end{punkt}

We also have so-called intermediate $S$-model structures: Let $S=\{A_\ell\to B_\ell|\ell\aus L\}$ be a set of cofibrations in \pre\, containing $I_{\rm proj}$, see \ref{gen sets for proj}. We denote by $I_S$ the set of cofibrations of the form
\begin{equation}\label{saturation}
      (A_\ell\times\Delta^n)\sqcup_{(A_{\ell}\times\partial\Delta^n)}(B_\ell\times\partial\Delta^n)\to B_\ell\times\Delta^n .
\end{equation}
The {\bf\boldmath$S$-cofibrations} ${\rm Cof}_S$ are given by the saturation of $I_S$. The {\bf\boldmath $S$-fibrations} are given by the right lifting property with respect to $S$-cofibrations that are also local weak equivalences. Observe:
    $$ {\rm Cof}_{\rm proj}\ent{\rm Cof}_S\ent{\rm Cof}_{\rm inj} $$
The fact, that this gives proper simplicial cofibrantly generated model structure on \pre\, is proved in \cite{jar:intermediate}.
We will refer to these model structures as {\bf intermediate model structures} including the projective and the injective ones and denote a choice of one of them by \mc{I}. 
\end{Def}

\begin{Bem} \label{all serve well}
The local or projective model structure on \pre\, with respect to the chaotic topology on \mc{C} was first constructed in \cite{BK:lim} and generalized to other topologies in \cite{Blander:local}. The global or injective structure was constructed in \cite{jar:simplicial-presheaves}. 
Its predecessor in the sheaf case was found in \cite{Joyal:letter}. 
The $S$-model structures on \pre\, were constructed in \cite{jar:intermediate}. Their existence can also be derived from a more general context in \cite{Beke:1}.

All of these structures serve equally well as a starting point for a theory of (pre-)\-sheaves of $n$-types in the next section.
\end{Bem}

Of course, these model structures are all Quillen equivalent to each other by the identity functor. To summarize we give the following theorem.
\begin{satz} 
Each \mc{I}-model structure from \emph{\ref{all serve well}} is a cofibrantly generated proper simplicial model structure on \pre.
\end{satz}

The interesting part here is cofibrant generation. In fact it makes everything work by the small object argument. It is easily seen for the projective structure and the generating sets are easy to display. First let us define the free presheaf functor $L_U$.

\begin{Def} \label{L_U}
For each $U\aus\mc{C}$ let $L_U\co\mc{S}\to\pre$ be the left Kan extension of the functor which assigns to every simplicial set the constant simplicial presheaf over $U$, i.e. for every $K\aus\mc{S}$ and every $V\aus\mc{C}$:
    $$ L_UK(V):=\coprod_{V\to U}K $$
So $L_U$ is the left adjoint to taking $U$-sections. 
\end{Def}

\begin{Bem} \label{gen sets for proj}
Now the generating sets for the projective model structures on \pre\, are given by:
\begin{center}
    $\begin{array}{cccl}
        I_{\rm proj} & := & \{\ \hfill L_U\partial\Delta^{s+1}\to L_U\Delta^{s+1}\ \hfill | & U\aus\mc{C}, s\bth 0\ \} \cup \{\varnothing \to L_U\Delta^0\,|\,U\aus\mc{C}\} \\
        J_{\rm proj} & := & \{ \hfill L_U\Lambda^{s+1}_k\to L_U\Delta^{s+1} \hfill | & U\aus\mc{C}, s\bth 0,\ s+1\bth k\bth 0 \} 
    \end{array}$
\end{center}
The projective model structure exists more generally for small presheaves on arbitrary sites as described in \cite{Chorny-Dwyer:small} and \cite{BCR:calc}. The real challenge lies in the injective structure. Then the case of the other intermediate model structures is done by juggling with what we already have \cite{jar:intermediate}.
\end{Bem}

\begin{Bem} \label{small site}
It is in the proof of the existence of the injective model structure on \pre, where we have to be careful with the set-theory. Let \mc{C} be a small site, i.e. with a set of objects. Let $\alpha$ be an infinite cardinal such that $\alpha\uber|{\rm Mor}\mc{C}|$. A presheaf $X$ is $\alpha$-bounded if
  $$ |\bigsqcup_{U\in\mc{C}}\bigsqcup_{n\bth 0}X_n(U)|\sth\alpha.$$
Define:
\begin{align*}
        I_{\rm inj} := & \{ i\co A\inj L_U\Delta^n\,|\,n\bth 0, i\hbox{ sectionwise monomorphism} \} \\
        J_{\rm inj} := & \{ j\co A\inj B\,|\, j\hbox{ local weak equivalence and sectionwise monomorphism}, \\ 
                       & B\,\, \alpha{\rm -bounded}\} 
\end{align*}
These sets are generators for the (trivial) cofibrations of the injective structure on \pre. The factorizations are constructed by the small object argument with respect to these sets. We will need these sets for the injective structure on \pregd\, and \shvgd.

The description of the generating set of trivial cofibrations for the intermediate $S$-structures is a little bit involved and to that end we refer to \cite{jar:intermediate}. Let us just note it exists and call it $J_S$. The set $I_S$ was described in \Ref{saturation}.
\end{Bem}

Note that all fibrations, no matter in which model structure on \pre, are in particular objectwise fibrations. Hence the following theory of long exact homotopy sequences, which we have copied from \cite{jar:lec-simp-pre}, is always available.
\begin{Bem} \label{long exact homotopy sequence}
A projective fibration $p\co X\to Y$ of projectively fibrant simplicial presheaves induces a long exact sequence of homotopy sheaves in the following way:
There is a simplicial presheaf $F\to X_0$ fibred over $X_0$ such that
    $$ F(U) = \bigsqcup_{x\in X_0(U)} p^{-1}p(x) $$
for all $U\aus\mc{C}$. Write $F_{p(x)}=p^{-1}p(x)$. The object $F$ has homotopy groups $\pi_sF\to X_0$ with
    $$ \pi_sF(U)=\bigsqcup_{x\in X_0(U)}\pi_s(F_{p(x)},x) $$
for $s\bth 1$, and there is a fibred sheaf of sets $\wt{\pi}_0F\to \wt{X}_0$ defined in the obvious way. There are also homotopy groups $\pi_sX\to X_0$ and $\pi_sY_p\to X_0$, where by the latter we mean
    $$ \pi_sY_p=\bigsqcup_{x\in X_0(U)}\pi_s(Y,p(x)), $$
which is fibred in an obvious way over $X_0$. For $s\bth 1$ these fibred group objects carry an action by the fundamental groupoid $\pi X$ of $X$, and there is a long exact sequence of $\pi X$-functors
    $$ \dots\to\pi_2Y_p\toh{\partial}\pi_1F\to\pi_1X\to\pi_1Y_p\to\dots $$
fibred over $X_0$. 
As depicted we will usually leave the map to $X_0$ to be understood.
We can extend this exact sequence to the right. Write $p^*\pi Y$ for the inverse image of $\pi Y$ along the groupoid map $\pi p\co\pi X\to\pi Y$. Then ${\rm Aut}_{p^*\pi Y}(x)=\pi_1(Y,p(x))$. The path components $\pi_0F_{p(x)}$ of the fibres determine a functor $\pi_0F\co p^*\pi Y\to {\rm Sets}$. This functor restricts to the standard action of $\pi_1(Y,p(x))$ on $\pi_0F_{p(x)}$, and also restricts to a functor $\pi_0F\co \pi X\to {\rm Sets}$. The usual boundary maps
    $$ \pi_1(Y,p(x))\to\pi_0F_{p(x)} $$
are $\pi X$-equivariant, meaning that they determine a map $\partial\co\pi_1Y\to\pi_0F$ of $\pi X$-funtors. The set $\pi_0F_{p(x)}$ is pointed be $[x]$, so that the map $\partial$ is a transformation of $\pi X$-functors taking values in pointed sets. One shows that the sequence
    $$ \pi_1X\toh{p_*}\pi_1Y\toh{\partial}\pi_0F$$
is exact in each section. The diagrams
\begin{center}
    $\xymatrix@=10pt{ \pi_0F_{p(x)} \ar[dd]_{\alpha_*}\ar[dr] & \\
                    &  \pi_0X \ar[dl] \\
                \pi_0F_{p(y)}  }$
\end{center}
commute for each $\alpha\co p(x)\to p(y)$ in $\pi Y$ and the induced function
    $$ \colim_{x\in p^*\pi Y} \pi_0F_{p(x)}\to\pi_0X $$
is a bijection. The ``action'' of the groupoid $p^*\pi Y$ on the various $\pi_0F_{p(x)}$ restricts to an action of $\pi_1(Y,p(x))$ on the set $\pi_0F_{p(x)}$, and there is a commutative diagram
\diagr{ \pi_0F_{p(x)} \ar[r]^-\pi\ar[d] & \pi_0F_{p(x)}/\pi_1(Y,p(x)) \ar[d]^i \\
        \lim_{y\in p^*\pi Y}\pi_0F_{p(y)} \ar[r]_-\cong & \pi_0X }
where $\pi$ is the canonical surjection and $i$ is an injection.
Note  that we can now sheafify everything to obtain such a long exact sequence for the associated homotopy sheaves. 
\end{Bem}

Let us finally describe the theory of simplicial sheaves. We denote by
    $$ i\co\shv\leftrightarrows\pre:\!L^2 $$
the adjoint pair given by the forgetful functor $i$ and the sheafification functor $L^2$.

For a class \mc{C} of morphisms in a category we denote by {\bf\boldmath\mc{C}-proj} the class of morphisms having the right lifting property with respect to all elements of \mc{C}, and by {\bf\boldmath\mc{C}-cof} the class of morphisms having the left lifting property with respect to all elements of \mc{C}-proj.

\begin{Def} \label{sheaf version}
A map in \shv\, is called a {\bf local weak equivalences} or an {\bf injective cofibration} if it is one in \pre. 
The maps of the class $L^2(I_{\rm proj})-\text{cof}$ are called {\bf projective cofibrations}.  
Given a set $S$ of cofibrations in \shv\, containing $L^2(I_{\rm proj})$ the class of {\bf\boldmath $S$-cofibrations} is given by the saturation of $I_S$ in \shv\, constructed as in \ref{Def. of cofibrations}.
The corresponding fibrations in all cases are defined by the right lifting property.
\end{Def}

\begin{satz} 
Each of the model structures from \emph{\ref{sheaf version}} is a cofibrantly generated proper simplicial model structure on \pre. The functors
    $$ i\co\shv\leftrightarrows\pre:\!L^2 $$
form a Quillen equivalence for the respective structures on both sides.
\end{satz}

The proof of the lifting axioms in this theorem relies on the useful observation proved in \cite{jar:lec-simp-pre}, that a map is a fibration in one of the model structures from \ref{sheaf version} in \shv\, if and only if it is a fibration in the corresponding structure on \pre. We will use the analogous argument in the proof of \ref{shvgd} for the groupoid case.

\section{Truncated simplicial presheaves and sheaves}
\label{section:Postnikov-Lokalisierung}

For ordinary simplicial sets there is a coaugmented functor $P_n\co\mc{S}\to\mc{S}$, called the $n$-th Postnikov section, whose coaugmentation induces the following isomorphisms
\begin{equation} \label{Postnikov stage}
      \pi_sP_nK\cong\left\{
                             \begin{array}{cl}
                                  \pi_sK &,\hbox{ for }0\sth s\sth n \\ 
                                     0   &,\hbox{ for }s\uber n \\ 
                             \end{array}
                     \right. 
\end{equation}
of homotopy groups for all basepoints of a Kan complex $K$.
Since the classical construction of $P_n$ is only homotopy meaningful on fibrant simplicial sets we precompose it with a fibrant replacement functor on \mc{S}, e.g. ${\rm Ex}^\infty$, and still denote the resulting functor by $P_n$. So for a simplicial set $K$, $P_nK$ denotes the $n$-th Postnikov section of ${\rm Ex}^\infty K$.
Now we are going to extend this to presheaves. Assume that we have chosen once and for all one of the model structures mentioned in remark \ref{all serve well}.

\begin{Def}
For a simplicial presheaf $X$ we can consider the presheaf $P_nX:=P_n\circ X\co\mc{C}^{\rm op}\to\mc{S}\toh{P_n}\mc{S}$. 
This induces a functor $P_n$ on simplicial presheaves which comes equipped with a natural transformation $p_n\co \id\to P_n$ from the identity functor to the $n$-th Postnikov stage, which for $n\bth 1$ induces isomorphisms
    $$ \wt{\pi}_sX\cong\wt{\pi}_sP_nX $$
for $0\sth s\sth n$ of sheaves fibred over $\wt{X}_0$ and isomorphisms
    $$ \wt{\pi}_sP_nX\cong\wt{X}_0 $$
for $s\bth n$. 
\end{Def}

Given the definition \ref{local weak equivalences} about the relation of the classical case with the presheaf case, these equations are the exact analogs of the equations \Ref{Postnikov stage}.
The coaugmented functor $P_n\co\pre\to\pre$ satisfies the axioms {\bf (A.4)}, {\bf (A.5)} and {\bf (A.6)} given in \ref{BF-axioms}.

\noindent
{\bf Proof of the axioms:} {\bf (A.4)} and {\bf (A.5)} are obvious. Axiom {\bf (A.6)} follows immediately from the long exact sequence of homotopy groups of the pullback square: 
     $$ ...\to \wt{\pi}_sA\to \wt{\pi}_sX_\alpha\oplus \wt{\pi}_sB_{g}\to \wt{\pi}_sY_{f\alpha}\to \wt{\pi}_{s-1}A\to ... $$
fibred over $A_0$. 
The isomorphisms $\wt{\pi}_sB\cong\wt{\pi}_sY_\beta$ for $0\sth s\sth n$ fibred over $B_0$ imply $\wt{\pi}_sA\cong \wt{\pi}_sX_\alpha$ in the same range fibred over $A_0$.
{\par \hfill $\Box$ \hspace{1cm} \par}

Then we obtain the following model structure.

\begin{Def} 
For any intermediate model structure \mc{I} on \pre\, the {\bf \boldmath$n$-\mc{I}-model structure} on \pre\, is given by the following classes of maps: 
A cofibration will simply be an \mc{I}-cofibration of the underlying intermediate structure. Further we will call a morphism $X\to Y$ in \pre
\begin{punkt}
    \item an {\bf \boldmath$n$\unboldmath-equivalence} if the induced map $P_nX\to P_nY$ is a local weak equivalence.
    \item an {\bf \boldmath$n$-\mc{I}\unboldmath-fibration} if it has the right lifting property with respect to all cofibrations of the underlying structure that are at the same time $n$-equivalences.  
\end{punkt}
We will write $\pre^n$ for the category of simplicial presheaves equipped with the $n$-\mc{I}-model structure.
\end{Def}

Obviously the $n$-equivalences are exactly the maps $f\co X\to Y$ that induce isomorphisms $\pi_0X\cong\pi_0Y$ and 
    $$ \wt{\pi}_sX = \coprod_{x\in X_0}\wt{\pi}_s(X,x)\to\coprod_{x\in X_0}\wt{\pi}_s(Y,fx)=\wt{\pi}_sY_f $$
of fibred objects for $1\sth s\sth n$.

\begin{satz} 
The $n$-\mc{I}-structure on \pre\, arising from any intermediate structure \mc{I} is a proper simplicial model structure. A map $p\co X\to Y$ is an $n$-\mc{I}-fibration if and only if it is an \mc{I}-fibration and the following diagram
\diag{ X \ar[r]^{p_n(X)}\ar[d]_p & P_nX \ar[d]^{P_np} \\
       Y \ar[r]_{p_n(Y)} & P_nY  }{P_n-pullback}
is a homotopy pullback square in the original structure.
\end{satz}

This theorem follows directly from \ref{Bousfield-Friedlander}, that is from \cite[9.3.]{Bou:telescopic}. The simplicial structure is the one of the underlying category described in \ref{simplicial structure}. If we decide to change the underlying model structure, then the $n$-model structure arising from the new intermediate structure will be Quillen equivalent to the old one by the identity functor.

It is worth mentioning that we obtain the same truncated model structures if we use the coskeleton functor instead of Postnikov sections: cofibrations are the same anyhow, and since $\cosk_{n+1}K\simeq P_nK$ for fibrant $K$ weak equivalences are the same. So in the pullback diagram \Ref{P_n-pullback} we can replace $P_n$ by $\cosk_{n+1}$.
The coskeleton is a purely categorical construction. However it does not preserve projective fibrations as $P_n$ does. Therefore we prefer the latter one.

We have the following characterization of $n$-\mc{I}-fibrations.

\begin{kor} \label{Post-n-fib}
A map $p\co X\to Y$ is an $n$-\mc{I}-fibration if and only if it is an \mc{I}-fibration and for some projective fibrant approximation $\wt{X}\to\wt{Y}$ the induced maps $\wt{\pi}_s\wt{X}\to\wt{\pi}_s\wt{Y}_p$ are isomorphisms of objects fibred over $\wt{X}_0$ for all $s\uber n$.
\end{kor}

\begin{beweis}
This follows from the comparison between the long exact homotopy sequences of the vertical maps in the homotopy pullback square \Ref{P_n-pullback}.
\end{beweis}

It follows that the \mc{I}-fibrant objects of the $n$-\mc{I}-structure are the \mc{I}-fibrant simplicial presheaves $X$ such that $\wt{\pi}_sX\cong\wt{X}_0$ for $s\uber n$. 

\begin{Bem} \label{local=sectionwise}
A map between globally fibrant objects is a local weak equivalence if and only if it is a sectionwise weak equivalence:
it is a weak equivalence between objects which are fibrant and cofibrant and is therefore a homotopy equivalence, and hence a weak equivalence in all sections.
This is proved explicitly in \cite{jar:lec-simp-pre}, although it was in the background of Jardine's papers for some time.
The following corollary is \cite[Prop. 6.11]{jardine:etale}, but we repeat the proof.
\end{Bem}

\begin{kor} 
A $n$-injectively fibrant object in \pre\, is an $n$-type in each section.
\end{kor}

\begin{beweis}
The proof is by induction. Suppose first that $F$ is $0$-injectively fibrant. So it is injectively fibrant and for $s\bth 1$ we have $\wt{\pi}_sF\cong\wt{F}_0$. This shows that $\wt{F}$ is isomorphic to $\wt{\pi}_0F$ viewed as a constant simplicial object. Then by the previous remark $F\simeq\wt{F}$ is sectionwise acyclic.

Now suppose $F$ is $1$-injectively fibrant, so that $\wt{\pi}_sF\cong\wt{F}_0$ for $s\bth 2$. Then $F$ is weakly equivalent to its first Postnikov section $P_1F$.  We have a fibre sequence
   $$ GK(\pi_1F,1)\to GP_1F \to GP_0F ,$$
where $G$ denotes a functorial injectively fibrant model. $GP_0F$ is sectionwise acyclic by the previous argument. Furthermore:
    $$ \pi_sGK(A,n)(U)\cong\begin{array}{cl}
                            H^{n-s}(U, A|_U) &,\hbox{ for } s\sth n \\
                                  0          &,\hbox{ for } s\uber n
                           \end{array} $$
So the long exact homotopy sequence proves that $F$ has homotopy groups sectionwise vanishing above degree $1$. We proceed inductively along the Postnikov tower.
\end{beweis}

\begin{Bem} \label{erzeugende Mengen}
Obviously, if $I_{\mc{I}}$ is a set of generating cofibrations for the underlying intermediate model structure \mc{I}, then $I_{\mc{I}}$ is also a set of generating cofibrations for the $n$-structure arising from it, because we did not change cofibrations. 

We wish to describe generating sets of trivial cofibrations relative to the underlying intermediate structure \mc{I}. 
Let $J_{\mc{I}}$ be a set of generating trivial cofibrations for \mc{I}. For each $U\aus\mc{C}$ let $L_U\co\mc{S}\to\pre$ be the functor defined in \ref{L_U}
Then we set:
    $$  J_{\mc{I},n} := J_{\mc{I}}\cup\{ \hfill L_U\partial\Delta^{s}\to L_U\Delta^{s} \, |\, U\aus\mc{C}, s\bth n+2\}\cup\{\ast\to L_U\partial\Delta^{n+2}\,|\,U\aus\mc{C}\} $$
The sets $J_\mc{I}$ were described in \ref{small site}. The next result gives $n$-structures a very combinatorial flavour. Interestingly enough, this result will not be used in the rest of the article except of providing similar generating sets for the truncated structures on groupoids. 
\end{Bem}

\begin{lemma} \label{J_n=erzeugende triviale Kofaserungen}
The set $J_{\mc{I},n}$ forms a set of generating trivial cofibrations for the $n$-structure on \pre.
\end{lemma}

\begin{beweis}
Fix an intermediate model structure \mc{I}.
We will prove that a map has the right lifting property with respect to $J_{\mc{I},n}$ if and only if it is an $n$-\mc{I}-fibration. Let $p\co X\to Y$ be an \mc{I}-fibration. Since all intermediate model structures on \pre\, are right proper it suffices by the following lemma \ref{Kan-right-proper} to restrict our attention to maps $p$ between injectively fibrant simplicial presheaves. 

We invoke \ref{local=sectionwise}.
It is then clear by \ref{Post-n-fib} that an $n$-\mc{I}-fibration between injectively fibrant objects has the right lifting property with respect to $J_{\mc{I},n}$. 

Conversely, let $p$ have the right lifting property with respect to $J_{\mc{I},n}$. Then $J_{\mc{I}}$ accounts for the fact that $p$ is an \mc{I}-fibration, $\{L_U\partial\Delta^{s+2}\to L_U\Delta^{s+2}\,|\,s\bth n, U\aus\mc{C}\}$ for $\wt{\pi}_sF\cong\wt{X}_0$ for $s\bth n+1$, where $F$ is the fiber of $p$ as in \ref{long exact homotopy sequence}, and the right lifting property with respect to $\{\ast\to L_U\partial\Delta^{n+2}|\,U\aus\mc{C}\}$ shows that $\wt{\pi}_{n+1}X\to\wt{\pi}_{n+1}Y_p$ is an epimorphism.
\end{beweis}

In the previous proof we have used the following lemma due to Kan and whose proof we will display for the benefit of the reader. 
\begin{lemma} \label{Kan-right-proper}
Let \mc{M} be a right proper model category. Consider the following diagram
\diagr{ A\ar[d]_i\ar[r] & X \ar[d]^p\ar[r]^{\simeq} & \wt{X} \ar[d]^{\wt{p}} \\
        B\ar[r]\ar@{.>}[ur]^-{h_1}\ar@{.>}[urr]_->>>>{h_2} & Y \ar[r]_{\simeq} & \wt{Y}  }
where $i$ is a cofibration and $p$ and $\wt{p}$ are fibrations. Then a lifting $h_1$ exists if and only if a lifting $h_2$ exists.
\end{lemma}

\begin{beweis}
Existence of $h_1$ clearly implies the existence of $h_2$. Now suppose $h_2$ is given. Construct a pullback:
\diagr{ P \ar[r]\ar[d] & \wt{X} \ar[d]^{\wt{p}} \\
        Y \ar[r]_{\simeq} & \wt{Y} \ar@{}[ul]|->>{\pullback} }
The canonical map $X\to P$ is a weak equivalence by right properness. Factor it into a trivial cofibration $j\co X\to E$ followed by a trivial fibration $E\to P$. $X$ is a retract of $E$ over $Y$ by the following diagram:
\diagr{ X \ar@{=}[r]\ar[d]_j & X \ar[d]^p \\
        E \ar[r]\ar[ur]_r & Y }
In the following diagram the lifts are easily seen to exists, $h_1$ is given by $rH$:
\diagr{ A \ar[d]\ar[r] & X \ar[r] & E \ar@{->>}[r]^\simeq\ar@/_8pt/[l]_r & P \ar[r]\ar[d] & \wt{X} \ar[d] \\
        B \ar[rrr]\ar[urrr]\ar[urr]^H &&&  Y \ar[r] & \wt{Y} \ar@{}[ul]|->>{\pullback} }
\end{beweis}

Lemma \ref{J_n=erzeugende triviale Kofaserungen} allows us to characterize fibrations by point set data. As an example we have the following description of $n$-fibrations of simplicial sets, which seems to be unknown even in this elementary case. Fibrant $n$-types, i.e. fibrant objects in the $n$-structure, can therefore be described in a quite nice combinatorial manner. The lemma also cleans up a loose end from \cite[3.9.]{Biedermann:truncated}, where $M_sK$ denoted the $s$-th matching set of the simplicial set $K$.

\begin{lemma} \label{Postnikov-Faserungen}
Let $K\to L$ be a fibration between fibrant simplicial sets and $n\bth 0$. This map induces isomorphisms on homotopy groups in degrees $s\uber n$ for all basepoints, i.e. it is an $n$-fibration, if and only if for $s\bth n+2$ the induced maps
    $$ K_s \to M_sK\times_{M_sL}L_s $$
and the map
    $$ M_{n+2}K\to M_{n+2}Y\times_{Y_0}X_0$$
are surjective. Here the last map is induced by some map $\ast\to\partial\Delta^{n+2}$.
\end{lemma}

\begin{beweis}
This follows directly from lemma \ref{J_n=erzeugende triviale Kofaserungen} in the case of $\mc{C}=*$ by adjointness.
\end{beweis}

We can describe an analogous $n$-\mc{I}-model structure for the category \shv\, of simplicial sheaves.

\begin{Def} 
The {\bf\boldmath$n$-\mc{I}\unboldmath-structure} on \shv\, is given by the following classes of maps: We call a morphism $X\to Y$ in \shv
\begin{punkt}
    \item an {\bf \boldmath$n$\unboldmath-equivalence} if $X\to Y$ is an $n$-equivalence in \pre.
    \item an {\bf \boldmath$n$-\mc{I}\unboldmath-fibration} if it has the right lifting property with respect to all \mc{I}-cofibrations in \shv\, that are at the same time $n$-equivalences.  
\end{punkt}
The cofibrations of this structure will be the usual \mc{I}-cofibrations of \shv. 
\end{Def}

Along the lines of \cite{jar:lec-simp-pre} the following theorem is straightforward to check.
\begin{satz} 
The $n$-\mc{I}-structure on \shv\, is a proper simplicial model structure. The adjoint pair of functors
    $$ i\co\shv\leftrightarrows\pre:\!L^2 $$
form a Quillen equivalence if we provide both categories with the $n$-\mc{I}-structure.
\end{satz}

I would like to summarize quickly, how these $n$-structures for $0\sth n\sth\infty$ fit together. Everything said here applies to sheaves and presheaves of simplicial sets as well as later to groupoids enriched in simplicial sets.

The identity functor on \pre\, maps $n$-equivalences to $(n-1)$-equi\-va\-len\-ces and preserves cofibrations. So $\id\co\pre^n\to\pre^{n-1}$ is a left Quillen functor and we obtain a Quillen pair:
       $$ \id\co\pre^n\leftrightarrows\pre^{n-1}:\!\id $$
The right derived $R(\id)$ is the composition $P_{n-1}$, where $P_{n-1}$ is the $(n-1)$-st Postnikov section composed with a fibrant approximation functor on \pre. Obviously, the functor
       $$ P_{n-1}\co\ho{\pre^{n-1}}\inj\ho{\pre^n} $$
as well as the functor  
       $$ P_{n}\co\ho{\pre^{n}}\inj\ho{\pre} $$
are embeddings of full subcategories. We can rephrase this by saying that $P_n$ is a coreflection.

\section{Sheaves and Presheaves of simplicial groupoids}
\label{section:groupoids}

Let Gd denote the category of small groupoids. Let $s{\rm Gd}$ denote the category of small groupoids enriched in simplicial sets. An object in this category is a simplicial groupoid, whose simplicial set of objects is discrete. By abuse of language we will often refer to these gadgets simply as simplicial groupoids.
For a small Grothendieck site \mc{C} let $s{\rm Pre}_{\mc{C}}({\rm Gd})$ denote the category of presheaves in $s{\rm Gd}$. We will sometimes refer to the objects of this category simply as presheaves of simplicial groupoids. Analogously we will talk about sheaves of simplicial groupoids and the category \shvgd. 

We would like to put a model structure on \pregd, which is Quillen equivalent to a chosen one on \pre. The natural way to do this, is to transfer the model structure on \pre\, via a pair of adjoint functors. There are at least two such pairs. The first one is the adjoint pair of functors given by
    $$ \Pi d^*\co\pre\leftrightarrows s{\rm Pre}_{\mc{C}}({\rm Gd}):\!dB(\frei).$$
Here $B$ denotes the the bisimplicial presheaf obtained by applying the classifying space functor $B$ to a simplicial groupoid and $d$ denotes the diagonal of the bisimplicial object. Its left adjoint is the composition of the left adjoint $d^*$ of the diagonal functor $d$ with the fundamental groupoid functor $\Pi$. 
The second pair is  
    $$ G\co\pre\leftrightarrows s{\rm Pre}_{\mc{C}}({\rm Gd}):\!\ol{W},$$
where $G\co\pre\to\pregd$ is the loop groupoid functor and $\ol{W}$ is the universal cocycle functor from \cite{DK:simp-groupoids}, also discussed in \cite[V.]{GoJar:simp}.
The first pair was used in \cite{joyal-tierney:shv-simp-groupoids}. The second pair was taken up by \cite{Luo:model-2-groupoids}.
I follow the second approach. 

Let \mc{I} denote an intermediate model structure on \pre\, as in \ref{all serve well}.
\begin{Def} \label{transferred I}
We call the following classes of maps the transfered \mc{I}-model structure or simply {\bf\boldmath \mc{I}-structure} on $s{\rm Pre}_{\mc{C}}({\rm Gd})$. A morphism $G\to H$ is 
\begin{punkt}
   \item a {\bf local equivalence} if $\ol{W}G\to\ol{W}H$ is a local weak equivalence in \pre.
   \item an {\bf\boldmath \mc{I}-fibration} if $\ol{W}G\to\ol{W}H$ is an \mc{I}-fibration in \pre.
   \item a {\bf\boldmath \mc{I}-cofibration} if it has the left lifting property with respect maps which are \mc{I}-fibrations and local weak equivalences. 
\end{punkt}
\end{Def}

Before we embark on a proof of the existence of the model structure on \pregd\, we would like to describe generating sets for the (trivial) cofibrations. 
Note that $\ol{W}\co\pregd\to\pre$ was defined sectionwise, therefore it is obvious that we obtain generating sets for the \mc{I}-structure on \pregd\, by applying $G\co\pre\to\pregd$ to the corresponding generating sets in \pre.

To summarize we define:
\begin{align*}
   I_{\rm proj}^{\pregd} = \{& GL_U\partial\Delta^n\to GL_U\Delta^n\,|\,U\aus\mc{C} \} \\
   J_{\rm proj}^{\pregd} = \{& GL_U\Lambda^n_k\to GL_U\Delta^n\,|\,U\aus\mc{C}, n\uber 0, n\bth k\bth 0 \} \\
   I_{\rm inj}^{\pregd}  =\{& Gi\co GA\inj GB\,|\,\, i\,\, \hbox{\rm monomorphisms in }\pre, B\,\,\alpha{\rm -bounded}\} \\
   J_{\rm inj}^{\pregd}  = \{& Gj\co GA\inj GL_U\Delta^n\,|\,\, j\,\, \hbox{\rm monomorphisms and local weak equivalence} \\
                             & n\uber 0\} 
\end{align*}   
These are the generating sets for the projective and the injective structure on \pregd.
Since I did not describe the generating (trivial) cofibrations for the intermediate $S$-structure on \pre\, earlier, I will not describe the corresponding ones for \pregd\, now. I refer to \cite{jar:intermediate} and apply $G$.

\begin{satz} \label{model structures on pregd}
For an intermediate model structure \mc{I} from \emph{\ref{all serve well}} the transferred \mc{I}-structure on $s{\rm Pre}_{\mc{C}}({\rm Gd})$ from \emph{\ref{transferred I}} is a right proper cofibrantly generated model structure. The pair
    $$ G\co\pre\leftrightarrows s{\rm Pre}_{\mc{C}}({\rm Gd}):\!\ol{W} $$
forms a Quillen equivalence. 
\end{satz}

Unfortunately I do not know whether these model structures are left proper or simplicial. \\

\begin{beweis}
To prove the existence for the injective model structure there is \cite{Luo:model-2-groupoids}. The path to success is the following: first one proves that trivial injective cofibrations in \pregd\, are closed under pushouts by a Boolean localization argument. Then one uses the small object argument to prove the factorizations. Then the lifting axioms are either obvious or proved by the retract argument.

This can be bootstrapped to the other cases: We need to supply generating sets of (trivial) cofibrations, that determine the corresponding (trivial) fibrations via the lifting property. We have that as we saw just before stating the theorem. Then we realize, that (trivial) cofibrations push out: First pushouts of cofibrations are cofibrations by their definition through the right lifting property. Then, if they are trivial, they are in particular trivial injective cofibrations, which push out to weak equivalences as pointed out above \cite{Luo:model-2-groupoids}. Together this implies that (trivial) cofibrations are closed under pushouts, so the small object argument proves everything.

Right properness is clear for the projective structure, since all objects are projectively fibrant. Since all \mc{I}-fibrations are projective fibrations, right properness now follows in general.
\end{beweis}

Let us put the intermediate model structures on \shvgd. 
\begin{Def} 
A map in \shvgd\, will be called a {\bf local weak equivalence} or a {\bf projective}, or {\bf injective}, or {\bf\boldmath$S$-fibration} if it is so in \pregd. The corresponding cofibrations are always defined by the left lifting property.  
\end{Def}

The following theorem is proved in the same way as for the case of presheaves, see \cite{jar:simplicial-presheaves} or  \cite{jar:lec-simp-pre}. 
We simply realize that sets of generating (trivial) cofibrations are obtained by applying the sheafification functor $L^2$.
The second Quillen equivalence involving $(G,\ol{W})$ follows from \cite{DK:simp-groupoids}.
 
\begin{satz} \label{shvgd}
The intermediate model structures on \shvgd\, exist. They are right proper and cofibrantly generated. The pairs of adjoint functors
    $$ i\co\shvgd\leftrightarrows\pregd:\!L^2$$
and
    $$ G\co\shv\leftrightarrows\shvgd:\!\ol{W}$$
form a Quillen equivalences.
\end{satz}

The last theorem puts the last corner of the diagram \Ref{main} for the case $n=\infty$ in place. The fact that the diagram \ref{main} is commutative is obvious. We can now start to truncate the groupoid side of it.

\section{Truncated sheaves and presheaves of simplicial groupoids}
\label{section:truncated groupoids}
 
We can localize these structures on \pregd\, in the same way as we did with \pre. Let us describe the Postnikov section functor for simplicial groupoids. 
\begin{Def} 
For a simplicial groupoid $G$ we define its {\bf\boldmath $n$-th Postnikov section $P_nG$} in the following way:
   $$ {\rm Ob}(P_nG) = {\rm Ob}(G) $$
Then for all $x,y\aus{\rm Ob}(G)$ we take as morphisms:
   $$ (P_nG)(x,y) = P_n(G(x,y)), $$
where $P_n$ on the right hand side denotes the $n$-Postnikov section of a simplicial set. The source and target maps of $G$ induce canonical source and target maps for $P_nG$.
\end{Def}

\begin{Bem} 
For a simplicial groupoid $G$ there is a pullback diagram
\diag{ G(x,y) \ar[r]\ar[d] & {\rm Mor}(G) \ar[d]^{(s,t)} \\
       \ast   \ar[r]_-{(x,y)} & {\rm Ob}(G)\times{\rm Ob}(G) }{morphism pullback}
where $s$ and $t$ are the source and target maps ${\rm Mor}(G)\to{\rm Ob}(G)$. 
Hence we observe that $P_nG$ is simply given by the data $P_n{\rm Mor}(G)\rightrightarrows{\rm Ob}(G)$. 
Both vertical maps in \Ref{morphism pullback} are clearly projective fibrations. Note that for this reason we do not need to apply a fibrant replacement to give $P_n$ homotopy meaning.
\end{Bem}

The following observation can be found in \cite{jar:class-gerbes}. For a simplicial groupoid $G$ and an object $y$ in it let $G/y$ denote the slice category of $G$ over $y$, whose objects are the morphisms $x\to y$ and whose morphisms are the obvious commuting triangles. This is an honest simplicial groupoid, whose simplicial set of objects is not discrete. Still, $dB(G/y)$ is contractible.
\begin{lemma} \label{groupoid-pullback}
For all objects $x$ and $y$ of a simplicial groupoid $G$ there is the following homotopy pullback square:
\diagr{ G(x,y) \ar[r]\ar[d] & dB(G/y) \ar[d] \\
        \ast   \ar[r]_x & dBG }
\end{lemma}

In other words, this lemma states that for all $s\bth 1$ we have
\begin{equation} \label{shift}
   \wt{\pi}_sdBG \cong \wt{\pi}_{s-1}G
\end{equation}
fibred over $G_0$.
The previous lemma \ref{groupoid-pullback} can be extended to presheaves sectionwise.
\begin{lemma} \label{truncation for groupoids}
For $n\bth 0$ the canonical map $G\to P_n G$ induces isomorphisms
    $$ \wt{\pi}_sdBP_n G \cong \left\{ \begin{array}{cl}
                                           \wt{\pi}_sdBG &\hbox{ for }0\sth s\sth n+1, \\
                                           \wt{dBG}_0=\wt{G}_0 &\hbox{ else. }
                                          \end{array}
                                  \right.$$
\end{lemma}

\begin{beweis}
We certainly have $\wt{\pi}_0dBG\cong\wt{\pi}_0dBP_nG$. For all other $s\uber 0$ the claim follows from the isomorphisms in \Ref{shift}.
\end{beweis}

Lemma \ref{truncation for groupoids} states for $n\bth 1$:
\begin{equation}\label{shift2}
   \ol{W}P_{n-1}G\simeq dBP_{n-1}G\simeq P_{n}BG\simeq P_n\ol{W}G 
\end{equation}
It follows easily that the coaugmented functor $P_{n}\co \pregd\to s{\rm Pre}_{\mc{C}}({\rm Gd})$ satisfies the axioms {\bf (A.4)}, {\bf (A.5)} and {\bf (A.6)}, so that we can localize the category \pregd\, along it. 

\begin{Def} 
Again let \mc{I} denote an intermediate model structure on \pre\, from \ref{all serve well} and consider the transfered \mc{I}-structure on \pregd\, \ref{transferred I}.
For $n\bth 0$ we call the following classes of maps the {\bf\boldmath $n$-\mc{I}-structure} on \pregd. A morphism $G\to H$ is 
\begin{punkt}
    \item an {\bf\boldmath $n$-equivalence} if $P_nG\to P_nH$ is a local equivalence in \pregd.
    \item an {\bf\boldmath $n$-\mc{I}-fibration} if it has the right lifting property with respect to all \mc{I}-cofibrations in \pregd\, that are also $n$-\mc{I}-equivalences. 
\end{punkt}
Cofibrations remain the same as for the \mc{I}-structure on \pregd.
\end{Def}

\begin{satz} \label{presheaf equivalence}
The $n$-\mc{I}-structure on \pregd\, is a right proper cofibrantly generated model structure. For $n\bth 1$ the pair of functors
   $$ G\co\pre^n\leftrightarrows s{\rm Pre}_{\mc{C}}({\rm Gd})^{n-1}:\!\ol{W}$$
forms a Quillen equivalence. 
\end{satz}

\begin{beweis}
The existence of the model structure follows from the general theorem \ref{Bousfield-Friedlander}, since $P_n$ satisfies {\bf (A.4)}, {\bf (A.5)} and {\bf (A.6)}: The first two axioms follow from the corresponding fact about Postnikov section of simplicial sets and {\bf (A.6)} follows directly from the isomorphisms \Ref{shift2}.
Right properness is inherited. 

It follows from \ref{truncation for groupoids} that $\ol{W}$ maps $(n-1)$-\mc{I}-equivalences of presheaves of simplicial groupoids to $n$-\mc{I}-equivalences of simplicial presheaves.
Given the characterization in \ref{Bousfield-Friedlander} of fibrations in the localized model structures it also follows that $\ol{W}$ maps $(n-1)$-\mc{I}-fibrations in \pregd\, to $n$-\mc{I}-fibrations in \pre. Hence we have a Quillen pair. Then it follows, that a map $GX\to H$ is an $(n-1)$-\mc{I}-equivalence in \pregd\, if and only if the map
    $$ \ol{W}P_{n-1}G X\to\ol{W}P_{n-1}H $$
is a weak equivalence in \pre. But again from \ref{truncation for groupoids}:
    $$ \ol{W}P_{n-1}G X\simeq P_n\ol{W}G X\simeq P_nX$$
and $\ol{W}P_{n-1}G\simeq P_n\ol{W}G$. So we have a Quillen equivalence.

The generating set of trivial cofibrations is given by $GJ_{\mc{I},n}$, where $J_{\mc{I},n}$ from \ref{erzeugende Mengen} generates the trivial $n$-\mc{I}-cofibrations in \pre. We need to show that a map $H\to K$ in \pregd\, is an $(n-1)$-\mc{I}-fibration if it has the left lifting property with respect to $GJ_{\mc{I},n}$. But then $\ol{W}H\to\ol{W}K$ has the left lifting property with respect to $J_{\mc{I},n}$. This, together with the isomorphisms \Ref{shift2}, shows that the diagram
\diagr{ H \ar[r]\ar[d] & P_{n-1}H \ar[d] \\
        K \ar[r] & P_{n-1}K }
is a homotopy pullback. Since $H\to K$ is an \mc{I}-fibration, theorem \ref{Bousfield-Friedlander} shows, that $H\to K$ is an $(n-1)$-\mc{I}-fibration.
\end{beweis}

To cover the sheaf case we simply observe that we can copy the previous argument line by line. 
\begin{satz} 
The $n$-\mc{I}-structure on \shvgd\, is a right proper cofibrantly generated model structure. For $n\bth 1$ the pairs of functors
   $$ L^2\co\pregd^{n-1}\leftrightarrows\shvgd^{n-1}:\!i$$
and
   $$ G\co\shv^n\leftrightarrows\shvgd^{n-1}:\!\ol{W}$$
form Quillen equivalences. 
\end{satz}

Let $j_n\co\Delta_n\inj\Delta$ be the inclusion of the category of finite ordinals $\sth n$ into the whole category $\Delta$. For a complete category \mc{C} let $s_n\mc{C}$ denote the category of $n$-truncated simplicial objects over \mc{C}, i.e. functors $\Delta_n^{\rm op}\to\mc{C}$. Let further $j_n^*\co s\mc{C}\to s_n\mc{C}$ denote the restriction along $j_n$ and let $r^n\co s_n\mc{C}\to s\mc{C}$ be its right adjoint. Obviously we have
   $$ r_nj_n^*\cong\cosk_n.$$
As in the case of \pre\, and using the isomorphism
    $$ (\cosk_nG)(x,y)\cong\cosk_n(G(x,y))\simeq P_{n-1}G(x,y) $$
for every simplicial groupoid $G$, $n\bth 1$ and all $x,y\aus{\rm Ob}(G)$ we can truncate any \mc{I}-model structure on \pregd\, with respect to $\cosk_n$ and obtain the same $n$-\mc{I}-model structure. It is then easy to set up a Quillen equivalence 
    $$ j_n^*\co\pregd^n\leftrightarrows s_n{\rm Pre}_{\mc{C}}({\rm Gd}):\!r_n $$with actual $n$-truncated objects.

\bibliographystyle{alpha}

\end{document}